\documentclass{amsart}
\usepackage{amssymb,diagrams,placeins,wrapfig}

\ifx\pdftexversion\undefined
\usepackage[dvips]{graphics,color}
\usepackage[dvips,colorlinks]{hyperref}
\else
\usepackage[pdftex]{graphics,color}
\usepackage[pdftex,colorlinks]{hyperref}
\fi

\hypersetup{pdftitle={Non-abelian genus two Veech surfaces},
pdfauthor={Sergey Vasilyev}}

\usepackage{epsfig}
\graphicspath{{pics/}}

\theoremstyle{plain}

\newtheorem{lem}{Lemma}[section]

\newtheorem{cor}{Corollary}[section]

\newtheorem{thm}{Theorem}[section]

\numberwithin{equation}{section}

\theoremstyle{definition}
\newtheorem{defn}{Definition}[section]

\theoremstyle{remark}
\newtheorem{notation}{Notation}[section]

\newcounter{lempart}[lem]

\newcommand{\lempart}{\refstepcounter{lempart}\noindent{\bf(\alph{lempart})} }

\newcommand{\Id}{\mathrm{Id}}
\newcommand{\ord}{\mathop{ord}}
\newcommand{\M}[1][2]{\mathcal M_{#1}}
\newcommand{\QM}[1][2]{\mathcal{Q}M_{#1}}
\newcommand{\QQM}[1][2]{\mathcal{Q}M^{-}_{#1}}
\newcommand{\OM}[1][2]{\Omega M_{#1}}
\newcommand{\R}{\mathbb R}
\newcommand{\Z}{\mathbb Z}

\newcommand{\C}{\mathbb C}
\newcommand{\CP}{\mathbb C\mathbb P}
\newcommand{\lto}{\longrightarrow}

\renewcommand{\phi}{\varphi}
\renewcommand{\epsilon}{\varepsilon}
\renewcommand{\emptyset}{\varnothing}

\begin{document}

\title[Genus two Veech surfaces]{Genus two Veech surfaces arising from general quadratic differentials}
\author{Sergey Vasilyev}
\address{Department of Mathematics\\ University of Chicago\\ 5734 S. University Ave\\ Chicago, Illinois 60615\\USA}
\email{vasilyev@math.uchicago.edu}

\begin{abstract}We study Veech surfaces of genus 2 arising from quadratic
	differentials that are not squares of abelian differentials. We prove
	that all such surfaces of type $(2,2)$ and $(2,1,1)$ are arithmetic. In
	$(1,1,1,1)$ case, we reduce the question to abelian differentials of
	type $(2,2)$ on hyperelliptic genus~3 surfaces with singularities at
	Weierstrass points, and we give an example of a non-arithmetic Veech
	surface.  \end{abstract}

\maketitle

\section{Introduction}\label{Intro}

A Veech surface is a Riemann surface with a quadratic differential, such that
the derivatives of its affine deformations with respect to the flat structure
induced by the quadratic differential, form a lattice in $PSL_2(\R)$.
Classification of Veech surfaces in genus 2 has been the subject of several
recent works by McMullen (\cite{McMSpin}, \cite{McMDecagon}) and Calta
(\cite{Calta}). These works classify all genus 2 Veech surfaces given by
quadratic differentials that are squares of abelian differentials. In this
paper, we attempt to classify genus 2 Veech surfaces given by a general
quadratic differential.

Given a quadratic differential on a Riemann surface one can construct
a double cover on which this quadratic differential pulls back to a square of
an abelian differential. This provides a way to reduce the study of quadratic
differentials to the study of abelian differentials on Riemann surfaces of
higher genus. We extend this construction, and use it to describe genus 2 Veech surfaces given by general quadratic differentials.

We prove that all Veech surfaces of genus 2 given by quadratic differentials
with two double zeroes that are not squares of abelian differentials arise from
tori:

\begin{thm} All Veech surfaces in $\QQM(2,2)$ are arithmetic.
\end{thm}

Similarly all Veech surfaces of genus 2 with one double zero and two simple zeroes arise from tori:

\begin{thm}All Veech surfaces in $\QM(2,1,1)$ are arithmetic.
\end{thm}	

In the case of Veech surfaces of genus 2 with four simple zeroes, we reduce the
question to abelian differentials on genus 3 Veech surfaces:

\begin{thm}
There is a one-to-one correspondence between Veech surfaces of genus 2 in $\QM(1,1,1,1)$ and hyperelliptic Veech surfaces of genus 3 in $\OM[3](2,2)$ with singularities at Weierstrass points.
\end{thm}

\medskip {\it Acknowledgments.} I wish to thank my adviser Howard Masur for
introducing me to this area and proposing the original problem. I am grateful
to him for the many hours that we spent discussing this and other related
topics, and for the numerous comments on the original version of this paper. I
also wish to thank Alex Eskin for reading and commenting on the draft, and I
would like to acknowledge him for suggesting an equivalent alternative to the
construction in \autoref{MainConstruction} which lead to \autoref{thm:Eskin}.

\section{Background}\label{Prelim}

\subsection{Translation structures}

A {\it translation structure} on a Riemann surface $X$ is an atlas of coordinate
charts $\{(U_i,\phi_i:U_i\to\C)\}$ covering $X$ except maybe for some finite
set of points $\{S_1,S_2,\ldots S_k\}$, such that all transition maps
$\phi_j\circ\phi_i^{-1}:\C\to\C$ are translations $z\mapsto z+a$. Weakening
this condition to allow all maps of the form $z\mapsto \pm z+a$, we obtain a
{\it half-translation structure} on $X$. Since these transition maps preserve
families of parallel lines in $\C$, we obtain foliations $F_\phi$ of $X$ for
every direction $\phi$. In case of a translation structure, these foliations
are orientable. These foliations will possibly have singularities at points
$S_1,S_2,\ldots,S_n$. By counting the number of prongs and multiplying it by
$\pi$ we can assign cone angles to singularities.

An important class of translation surfaces arises from polygonal billiard
tables in which each angle is a rational multiple of $\pi$. Such a billiard
table defines a translation surface via an unfolding construction~(\cite{MT}).

\subsection{Quadratic differentials}

There is another way to describe (half-)translation structures. In case of a
translation structure, differentials $dz$ in each chart paste together to give
a holomorphic differential $\omega$ on $X\backslash\{S_1,S_2,\ldots,S_n\}$.
This differential can be extended to $X$: a cone singularity of angle $2k\pi$
will give rise to a zero of $\omega$ of order $k-1$. In case of a
half-translation structure, quadratic differentials $dz^2$ in each chart paste
together to give a (possibly meromorphic) quadratic differential $q$ on $X$
with zeroes of order $c-2$ at each cone singularity of angle $c\pi$. If $c=1$
then $q$ will have a pole of order $1$. 

Conversely, an abelian differential $\omega$ on $X$ defines a translation
structure by considering charts in which $\omega$ is given by $dz$.  Similarly,
a meromorphic quadratic differential with poles of order not larger than $1$
defines a half-translation structure on $X$. We would like to point out that in
this paper, a quadratic differential will mean a holomorphic quadratic
differential, unless it is specified otherwise. 

Denote by $\OM[g]$ (resp. $\QM[g]$) the moduli space of genus $g$ Riemann
surfaces with a choice of an abelian (resp. quadratic) differential. These
moduli spaces are further stratified by the orders of zeroes of the
corresponding differentials. These strata will be denoted by
$\OM[g](\epsilon_1,\ldots,\epsilon_n)$ and
$\QM[g](\epsilon_1,\ldots,\epsilon_n)$, where $\epsilon_1,\ldots,\epsilon_n$
are the orders of zeroes. We will think of $\OM[g]$ as a sub-space of $\QM[g]$
via $(X,\omega)\mapsto(X,\omega^2)$. We will also use the notation
$\QM[g]^-=\QM[g]\backslash\OM[g]$ for the moduli space of quadratic
differentials that are not squares of abelian differentials.

Every element of $(X,q)\in\QM[g]$ can be thought of as a Riemann surface $X$
with a half-translation structure $(U_i,\phi_i)$ with no cone singularities of
angle $\pi$. An element $A\in PSL_2(\R)$ acts on $(X,q)$, by changing each
coordinate map $\phi_i$ to $A\circ\phi_i$.  This action preserves $\OM[g]$ and
the stratifications by the orders of zeroes.

Quadratic differentials on a Riemann surface can be naturally thought of as
elements of the co-tangent space to the surface in the moduli space $\M[g]$ of
Riemann surfaces of genus $g$. Using the Teichm\"uller metric on $\M[g]$
co-tangent space is identified with the tangent space. This way a quadratic
differential gives rise to a tangent vector to the surface in the moduli space.
The projection of $PSL_2(\R)$ orbit of $(X,q)\in\QM[g]$ to $\M[g]$ is precisely
the complex geodesic through $X$ with respect to the Teichm\"uller metric in the
direction given by $q$.

\subsection{Veech surfaces}

An affine group $Aff^+(X,q)$ of $(X,q)\in\QM[g]$ is the group of all
orientation-preserving diffeomorphisms $X\to X$ that are given by affine maps
in each chart of the half-translation structure. The linear parts of these
affine maps are the same up to multiplication by ${\pm Id}$. Moreover, their
determinant is~$1$ because the total surface area (in the metric defined by
$q$) is preserved. Hence we get a well defined map $D:Aff^+(X,q)\to PSL_2(\R)$.
The image of this map is denoted by $SL(X,q)$ and is called the {\it Veech
group} of $(X,q)$ (strictly speaking it should be denoted by $PSL(X,q)$, but we
will allow ourselves an abuse of notation here).  We have an exact sequence $$
0\lto Aut(X,q)\lto Aff^+(X,q)\stackrel{D}{\lto}SL(X,q)\to 0, $$ where
$Aut(X,q)$ is the group of all holomorphic automorphisms of $X$ preserving
quadratic differential $q$.  

In case in which $q$ is a square of an abelian differential $\omega$, the map
$D$ is a well-defined map to $SL_2(\R)$ and the group $SL(X,\omega)=SL(X,q)$ is
a subgroup of $SL_2(\R)$.

$(X,q)$ is called a {\it Veech surface} if $SL(X,q)$ is a lattice in
$PSL_2(\R)$, that is it is a discrete subgroup such that the quotient
$PSL_2(\R)/SL(X,q)$ has finite (hyperbolic) volume. $(X,q)$ is a Veech surface
if and only if its $PSL_2(\R)$ orbit in $\QM[g]$ is closed (this was proved by
Smillie, see \cite{Smillie} for a sketch of the proof).  Veech surfaces satisfy
the so called Veech dichotomy: geodesic flow in every direction is either
periodic or uniquely ergodic (\cite{VeechDichotomy}). 

A Veech surface $(X,q)$ is called {\it primitive} if it cannot be realized as a
branched cover over $(X',q')$, where the genus of $X'$ is lower than the genus
of $X$ and the quadratic differential $q'$ pulls back to $q$ under the covering
map. If $(X,q)$ is a Veech surface that is not primitive, then the
corresponding surface $(X',q')$ is also Veech (see \autoref{GJ}), and hence
$(X,q)$ can be constructed as a branched cover of a Veech surface of lower
genus. 

A Veech surface is called {\it arithmetic} if its Veech group $SL(X,q)$ is
commensurable to $SL_2(\Z)$. As was proved by Gutkin and Judge,
$(X,\omega)\in\OM[g]$ is arithmetic if and only if it is a branched cover of a
torus (\cite[Theorem 5.5]{GJ}). Similarly, $(X,q)\in\QM[g]$ is arithmetic if
and only if the double cover given by $q$ (see \autoref{MainConstruction}) is a
branched cover of a torus.

\subsection{Genus 2 Veech surfaces}\label{IntroVeech2}

\begin{figure}[ht]
\begin{center}
$$\includegraphics{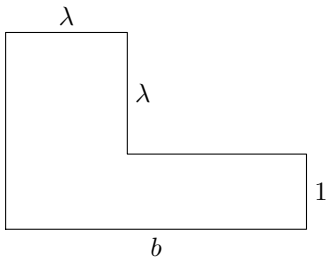}$$
\caption{An L-shaped billiard table determines a Veech surface, provided that $b\in\Z$, $\lambda=(e+\sqrt{e^2+4b})/2$, $e=-1,0$ or $1$, $e+1<b$, and if $e=1$ then $b$ is even.}
\label{fig:LShapedTable}
\end{center}
\end{figure}

Recent works by McMullen (\cite{McMSpin}) and Calta (\cite{Calta}) classify all
primitive Veech surfaces of genus~2 arising from abelian differentials with
a double zero. Up to the action of $SL_2(\R)$ all such surfaces are obtained
from an $L$-shaped billiard table of specific dimensions as illustrated in
\autoref{fig:LShapedTable} (\cite{McMSpin}).

Due to McMullen (\cite{McMDecagon}), in $\OM(1,1)$ there is a unique primitive
Veech surface up to the action of $SL_2(\R)$. It is obtained by gluing the
opposite sides of the regular decagon. McMullen's argument relies on a result
about torsion divisors proved by Moller (\cite[Cor. 3.4]{Mo}).

\subsection{Hyperelliptic surfaces}\label{sec:HyperellipticSurfaces}

One of the main tools in our study of genus 2 Veech surfaces is the
hyperelliptic involution. We would like to remind the reader of the main facts
about hyperelliptic surfaces that we will need later in the paper (for proofs
see \cite[Section III.7]{FarkasKra}).

A Riemann genus $g$ surface $X$ is hyperelliptic if there exists a two sheeted
covering $h:X\to\CP^1$. The corresponding sheet interchanging involution
$i_h:X\to X$ is called the hyperelliptic involution. If genus of $X$ is at
least 2, then the hyperelliptic involution is uniquely defined and does not
depend on the choice of $f$. Every surface of genus~$\leq 2$ is hyperelliptic.
For hyperelliptic surfaces, the Weierstrass points are the points fixed under
the hyperelliptic involution.  

One can obtain every hyperelliptic surface $X$ by starting with an affine curve
$$w^2=(z-z_1)(z-z_2)\ldots(z-z_{2g+2})$$ where all $z_j$'s are distinct, and
adding two points lying over $z=\infty$ by considering a second chart
$$z'=\frac1z\quad\mbox{and}\quad w'=\frac{w}{z^{g+1}}.$$ The hyperelliptic
involution is given by $(z,w)\mapsto (z,-w)$.  The Weierstrass points are
$(z_1,0),\ldots,(z_{2g+2},0)$.

In coordinates $(z,w)$ every holomorphic abelian differential
$\omega\in\Omega(X)$ is given by $$\frac{P(z)\,dz}w,\quad P(z)\in\C[z],\
deg(P)\leq g-1.$$  From this it is evident that $i_h(\omega)=-\omega$ for all
$\omega\in\Omega(X)$. It can also be shown that the zero divisor of
$\omega\in\Omega(X)$ is fixed by the hyperelliptic involution and the order of
zero at each Weierstrass point is even.  Conversely, any such divisor of degree
$g$ is a zero divisor of some holomorphic abelian differential.

The products of the holomorphic abelian differentials (taken two at a time)
form a $(2g-1)$-dimensional subspace of the $(3g-3)$-dimensional space of
quadratic differentials (\cite[p. 104, corollary 2]{FarkasKra}). Since $2g-1=3g-3$ for $g=2$, every quadratic
differential on a genus 2 surface can be written as a product of two abelian
differentials. In particular, every quadratic differential on a genus 2 surface
is fixed by the hyperelliptic involution.

\subsection{Strata of quadratic differentials in genus 2}

Suppose $(X,q)$ is a genus 2 surface with a quadratic differential $q$ that is
not a square of an abelian differential. The sum of orders of zeroes of $q$ on
$X$ is~4. Therefore there are several possible cases for zero configuration:
$(4)$, $(3,1)$, $(2,2)$, $(2,1,1)$ and $(1,1,1,1)$.

A single zero of order 4 is only possible if $q$ is a square of an abelian
differential. Indeed, $q$ is fixed by the hyperelliptic involution, and hence
the only zero of $q$ will have to be a Weierstrass point. Therefore there
is an abelian differential with a double zero at this point. The square of this
abelian differential will give a quadratic differential with the same zero
divisor as $q$, which means that it is proportional to~$q$. Therefore
$\QM(4)=\OM(2)$.

The zero configuration $(3,1)$ is impossible. If $q$ was such a quadratic
differential, then the hyperelliptic involution would have to fix the zero of
order 1, and would permute the three horizontal directions coming out of this
zero.  A permutation of order 2 on three elements has to fix one of the
elements, so the hyperelliptic involution would have to fix one of the
horizontal directions and would have to have infinitely many fixed points,
which is impossible. Therefore $\QM(3,1)=\emptyset$.

The goal of this paper is to study Veech surfaces in the strata $\QQM(2,2)$, $\QM(2,1,1)$ and $\QM(1,1,1,1)$.

\section{Reduction to abelian differentials}

\begin{notation} If $p:R\to S$ is a (ramified) double covering of Riemann
surfaces, the corresponding sheet-interchanging involution on $R$ will be
denoted by $i_p$.
\end{notation}

\subsection{Main Construction}\label{MainConstruction}

Let $X$ be a hyperelliptic surface and $q$ be a quadratic meromorphic
differential on $X$ with poles of order at most 1, which is not a square of an
abelian differential. Assume furthermore that $q$ is fixed by the hyperelliptic
involution (this is automatically satisfied if $X$ is a genus 2 surface and $q$
is a holomorphic quadratic differential; see \autoref{sec:HyperellipticSurfaces}). 

Let $f: Y\to X$ be the double covering given by $q$, i.e. $f^*(q)$ is a square
of an abelian differential $\alpha\in\Omega(Y)$ (cf. \cite[p.  519]{Lanneau}).
The covering map $f$ is branched over the odd-order zeroes and simple poles of
$q$. It follows from \autoref{lem:ZeroBehavior} below that $\alpha$ is a
holomorphic abelian differential even if $q$ has some simple poles. The
corresponding sheet-interchanging involution $i_f:Y\to Y$ sends $\alpha$ to
$-\alpha$.

Since $X$ is hyperelliptic, we have a ramified double covering $h_X:X\to\CP^1$.
The involution $i_{h_X}:X\to X$ is the hyperelliptic involution.

One can think of $Y$ as a set of pairs $(x,\tilde\alpha)$, where $x$ is a point
of $X$ and $\tilde\alpha$ is a holomorphic form defined locally around $x$,
s.t.  $\tilde\alpha^2=q$. The hyperelliptic involution $i_{h_X}$ acts naturally
on such pairs by sending $(x,\tilde\alpha)$ to
$(i_{h_X}(x),{i_{h_X}}_*(\tilde\alpha))$.  Indeed, by our assumption the
hyperelliptic involution $i_{h_X}$ preserves $q$, therefore it maps
$\tilde\alpha$ to another local square root of $q$.  Hence the hyperelliptic
involution on $X$ can be naturally lifted via $f$ to $Y$ to give an involution
$i_g$ on $Y$. Factoring $Y$ by this involution, we obtain a (ramified) double
covering $g:Y\to Z$.

\begin{lem}\label{ZHyperelliptic} $Z$ is hyperelliptic and $i_f$ descends
to $i_{h_Z}$, the hyperelliptic involution on $Z$ (i.e. $g\circ
i_f=i_{h_Z}\circ g$).
\end{lem}

\begin{proof} 
	
Pick a point $z\in Z$. Using notations above, there are two preimages of $z$
under $g$:  $(x,\tilde\alpha)$ and $(i_{h_X}(x),{i_{h_X}}_*(\tilde\alpha))$ for
some $x\in X$. Therefore $f\circ g^{-1}(z)$ consists of two (possibly
coinciding) points $x$ and $i_{h_X}(x)$, which are sent to the same point under
$h_X$. Therefore there is a well-defined map $h_Z:Z\to\CP^1$ making the
following diagram commute:

\begin{equation}
\label{diag:XYZ}
\begin{diagram}
&&Y&&\\
&\ldTo^f&&\rdTo^g&\\
X&&&&Z\\
&\rdTo_{h_X}&&\ldTo_{h_Z}&\\
&&\CP^1&&
\end{diagram}
\end{equation}

It remains to check that $h_Z$ is a (ramified) double covering. Take any
point $p\in\CP^1$, s.t. $h_X$ is not ramified at $p$ and $p$ is not an
image of a zero of $q$ under $h_X$. Then ${(h_X\circ f)}^{-1}(p)$ consists of four
points: $(x,\tilde\alpha)$,$(x,-\tilde\alpha)$,$(i_{h_X}(x),{i_{h_X}}_*(\tilde\alpha))$
and $(i_{h_X}(x),-{i_{h_X}}_*(\tilde\alpha))$:

\newarrow{Inv}<--->
\begin{equation}\label{diag:involutions}
\begin{diagram}
(x,-\tilde\alpha)&\rInv^{i_g}&(i_{h_X}(x),-{i_{h_X}}_*(\tilde\alpha))\\
\dInv <{i_f}&&\dInv >{i_f}\\
(x,\tilde\alpha)&\rInv^{i_g}&(i_{h_X}(x),{i_{h_X}}_*(\tilde\alpha))\\
\end{diagram}
\end{equation}

As indicated on the diagram (\ref{diag:involutions}), involution $i_g$
interchanges the columns while $i_f$ interchanges the rows. Therefore
${h_Z}^{-1}(p)=(g\circ{f}^{-1}\circ{h_X}^{-1})(p)$ consists of two points.
Moreover, $i_f$ descends to $i_{h_Z}$ on $Z$. Thus $h_Z$ defines a ramified
double covering of $\CP^1$, which proves that $Z$ is hyperelliptic.

\end{proof}

Consider involutions $i_f^*$ and $i_g^*$ on the space of holomorphic forms
$\Omega(Y)$.  As it is evident from the diagram (\ref{diag:involutions}),
involutions $i_f$ and $i_g$ commute on $Z$, and therefore $i_f^*$ and $i_g^*$
are commutative linear involutions on the linear space $\Omega(Y)$.  Hence
$\Omega(Y)$ decomposes into the sum of 1-dimensional subspaces, on which
$i_f^*$ and $i_g^*$ act by multiplication by $1$ or $-1$. If $i_f^*$ fixes a
form $\omega\in\Omega(Y)$, then $\omega$ descends to a holomorphic form on $X$.
But $i_g$ descends to the hyperelliptic involution on $X$.  Therefore
$i_g^*(\omega)=-\omega$. Similarly, if $i_g^*(\omega)=\omega$, then
$i_f^*(\omega)=-\omega$.  This discussion leads to the following conclusion:

\begin{lem}\label{OmegaYStructure} \lempart $\Omega(Y)=f^*(\Omega(X))\oplus g^*(\Omega(Z))$.\\
\lempart $i_f^*$ restricts to $\Id$ on $f^*(\Omega(X))$ and to $-\Id$ on $g^*(\Omega(Z))$.\\
\lempart $i_g^*$ restricts to $-\Id$ on $f^*(\Omega(X))$ and to $\Id$ on $g^*(\Omega(Z))$.
\end{lem}

\begin{cor}\label{GenusY} Genus of $Y$ is the sum of genera of $X$ and $Z$.
\end{cor}
\begin{proof} Follows from \autoref{OmegaYStructure}.
\end{proof}

\begin{cor}\label{cor:YHyperelliptic} $Y$ is a hyperelliptic surface with the
	hyperelliptic involution $i_{h_Y}$ given by $i_f\circ i_g$. In
	particular, the hyperelliptic involution of $Y$ descends to the
	hyperelliptic involutions on $X$ and $Z$.  \end{cor}

\begin{proof} By the lemma $(i_f\circ i_g)^*(\omega)=-\omega$ for any
holomorphic form $\omega$ on $Y$.  The first statement will follow from the following lemma.

\begin{lem} Assume $S$ is a Riemann surface and $i:S\to S$ is an involution,
such that $i^*(\omega)=-\omega$ for any holomorphic form $\omega$ on $S$.
Then $S$ is hyperelliptic, and $i$ is a hyperelliptic involution.
\end{lem}

\begin{proof} Consider $R=S/i$. If $R$ is not $\CP^1$ then there exists a non-zero holomorphic form $\omega$ on $R$.
But then the lift of $\omega$ to $S$ will be fixed by $i$, contradicting the hypothesis.
\end{proof}

Since $i_f$ is the sheet-interchanging involution of the covering $f:Y\to X$ and $i_g$ descends to $i_{h_X}$ on $X$, the hyperelliptic involution $i_{h_Y}=i_f\circ i_g$ descends to the hyperelliptic involution $i_{h_X}$. The same argument shows that $i_{h_Y}$ descends to $i_{h_Z}$ via $g:Y\to Z$.
\end{proof}

Since $i_f(\alpha)=-\alpha$, \autoref{OmegaYStructure} implies that $\alpha\in
g^*(\Omega(Z))$.  Therefore $\alpha$ descends to an abelian differential $\omega$ on
$Z$. This fact is fundamental to the rest of the paper. It has been known
before that one can reduce the study of a quadratic differential to a study of
an abelian differential on the corresponding double cover. Now, in case of a
hyperelliptic surface, we will be able to reduce this even further to a study
of an abelian differential on a surface of genus lower than the genus of the
corresponding double cover. We summarize the results of this section in the
following theorem:

\begin{thm}[Main Construction] 
	
Let $X$ be a hyperelliptic surface with a quadratic meromorphic differential
$q$ with poles of odd at most 1. Assume $q$ is not a square of an abelian
differential. Assume furthermore that $q$ is fixed by the hyperelliptic
involution of $X$ (this is automatically satisfied if genus of $X$ is 2 and $q$
is holomorphic).  Consider the double covering $f:Y\to X$ given by $q$.  The
quadratic differential $q$ lifts to a square of an abelian differential
$\alpha$ on Y.  The hyperelliptic involution of $X$ can be naturally lifted to
an involution $i_g$ of $Y$. Let $Z$ be the factor of $Y$ by the involution
$i_g$. Then $Z$ is a hyperelliptic Riemann surface, and $\alpha$ descends to an
abelian differential $\omega$ on $Z$.

\end{thm}

One can also obtain the same translation surface $(Z,\omega)$ from $(X,q)$ by
following the path $X\to\CP^1\gets Z$ in the diagram
(\ref{diag:XYZ}) on page \pageref{diag:XYZ}:

\begin{thm}\label{thm:Eskin} 
	
The quadratic differential $q$ descends via the map $h_X$ to a (possibly
meromorphic) quadratic differential $\check q$ on $\CP^1$. Then $(Z,\omega)$
can be obtained  from the quadratic differential $\check q$ on $\CP^1$ via the
double cover construction, i.e. $h_Z^*(\check q)=\omega^2$.

\end{thm}

\begin{proof}

The quadratic differential $q$ is fixed by the hyperelliptic involution $h_X$,
hence it descends to a quadratic differential $\check q$ on $\CP^1$. Using the
commutative diagram (\ref{diag:XYZ}) on page \pageref{diag:XYZ}, we see that
$g^*(h_Z^*(\check q))=f^*(h_X^*(\check q))=f^*(q)=\alpha^2=g^*(\omega^2)$.
Since $g:Y\to Z$ is a covering map, the map $g^*$ is injective on the spaces of
quadratic differentials. Therefore $h_Z^*(\check q)=\omega^2$.

\end{proof}

\subsection{Relation between \texorpdfstring{$(X,q)$ and $(Z,\omega)$}{(X,q) and (Z,w)}}

To establish a relation between the Veech groups of $(X,q)$ and $(Z,\omega)$,
we need to relate the Veech groups of the flat surfaces one of which covers the
other. This has been done independently by Gutkin and Judge (\cite{GJ}), and by Vorobets (\cite{Vorobets}).

\begin{defn} Let $(R,\tau)$ be a flat surface given by a quadratic differential
	$\tau$ on a Riemann surface $R$. Let $B$ be a finite subset of $R$.
	Denote by $(R,\tau,B)$ a flat surface obtained from $(R,\tau)$ by
	considering points in $B$ to be additional singularities. Consequently,
	$Aff^+(R,\tau,B)$ is the group of all orientation-preserving affine
	diffeomorphisms of $(R,\tau)$ mapping the set $B\cup Z(\tau)$ into
	itself, where $Z(\tau)$ is the set of zeroes of $\tau$, and
	$SL(R,\tau,B)$ is the group of derivatives of the diffeomorphisms in
	$Aff^+(R,\tau,B)$.

\end{defn}

\begin{thm}\label{GJ} Let $p:R\to S$ be
a ramified finite covering of Riemann surfaces branched over a finite set $B$. If
$\tau$ is a quadratic differential on $S$, then $SL(S,\tau,B)$ and
$SL(R,p^*(\tau))$ are commensurate.  \end{thm}

\begin{proof}

It was proved independently by Gutkin and Judge (\cite[Theorem 4.9]{GJ}), and
by Vorobets (\cite[Theorem 5.4]{Vorobets}) that in the case of coverings
ramified at singularities, the Veech groups of the base and the cover are
commensurate.  Therefore, $SL(S,\tau,B)$ and $SL(R,p^*(\tau),p^{-1}(B))$ are
commensurate (cf. \cite[Lemma 3]{HS}). It remains to be noticed that
$p^*(\tau)$ vanishes at the preimage of each branching point, and therefore
$SL(R,p^*(\tau),p^{-1}(B))=SL(R,p^*(\tau))$.

\end{proof}

\autoref{GJ} implies that $SL(X,q)$ and $SL(Y,\alpha)$ are commensurate, and
hence $(X,q)$ is a Veech surface if and only if $(Y,\alpha)$ is. Indeed, the
set $B_f$ of branching points of $f:Y\to X$ is the subset of zeroes of $q$, and
therefore $SL(X,q,B_f)=SL(X,q)$.  However, the set $B_g$ of
branching points of $g:Y\to Z$ is not necessarily a subset of zeroes of
$\omega$.  Applying \autoref{GJ} to $g:Y\to Z$, we obtain the main tool for
finding hyperelliptic Veech surfaces given by quadratic differentials:

\begin{thm} \label{MainThm}
	
Let $B_g$ be the set of branching points of covering map $g:Y\to Z$. Then
$SL(X,q)$ is commensurate to $SL(Z,\omega,B_g)$. Hence $(X,q)$ is a
Veech surface if and only if $(Z,\omega)$ is a Veech surface and
$SL(Z,\omega,B_g)$ is a finite index subgroup of $SL(Z,\omega)$.

\end{thm}

We will establish in the subsequent sections when $SL(Z,\omega,B_g)$ is a
finite index subgroup of $SL(Z,\omega)$ (see \autoref{lem:SL(Torus,B)} and
\autoref{lem:SL(H,B)}.)

\subsection{Reconstructing \texorpdfstring{$(X,q)$ from $(Z,\omega)$}{(X,q) from(Z,w)}}
\label{MainReconstruction}

We will need to reverse the main construction. For this take any holomorphic
differential $\bar\omega$ on $X$. We can lift it to a holomorphic differential
$\bar\alpha$ on $Y$. By \autoref{OmegaYStructure},
$i_g(\bar\alpha)=-\bar\alpha$. Therefore $\bar\alpha^2$ descends to a possibly
meromorphic quadratic differential $\bar q$ on $Z$ with poles of order at most
1 (see \autoref{lem:ZeroBehavior}). Since $\bar\alpha$ is fixed by $i_f$ and
$i_f$ descends to the hyperelliptic involution of $Z$
(\autoref{ZHyperelliptic}), we get that $\bar q$ is fixed by the hyperelliptic
involution of $Z$. Hence we can apply the main construction to $(Z,\bar q)$.
Since $i_f$ descends to the hyperelliptic involution on $Z$, the lift of the
hyperelliptic involution of $Z$ to $Y$ as described in
\autoref{MainConstruction} is either $i_f$ or $i_f\circ i_g$. We know that this
lift has to preserve $\bar\alpha$, hence it has to be $i_f$. Therefore the main
construction applied to $(Y,\bar q)$ will produce $(X,\bar\omega)$.

In other words, we have shown that there is a one-to-one correspondence between
triples $(X,q,\bar\omega)$ and $(Z,\omega,\bar q)$, where $q$ and $\bar q$ are
fixed by hyperelliptic involutions of $X$ and $Z$ correspondingly and have
poles of order at most 1.

\subsection{Some technical results}

To apply \autoref{MainThm} we need a description of the branching points
of $g:Y\to Z$ or, equivalently, the fixed points of $i_g$.

\begin{lem}\label{IgFixedPoints} The set of fixed points of $i_g$ consists of pre-images under $f$ of
all zeroes of $q$ of order not divisible by 4 that are also Weierstrass points.
\end{lem}
\begin{proof} 
	
Recall that we can think of $Y$ as a set of pairs $(x,\tilde\alpha)$, where
$\tilde\alpha$ is one of the two locally defined around $x\in X$ square roots
of $q$. Then $i_g$ sends $(x,\tilde\alpha)$ to
$(i_{h_X}(x),{i_{h_X}}_*(\tilde\alpha))$. For this point to be fixed under
$i_g$, we first of all need that $i_{h_X}(x)=x$, which means that $x$ is a
Weierstrass point. The map $i_g$ permutes the fiber over $x$. If $x$ is a
branching point, i.e. $x$ is an odd-order zero of $f$, then the fiber over $x$
consists of one point, and hence this point is fixed under $i_g$.  Otherwise
choose a local complex coordinate $z$ on $X$, s.t.  $z(x)=0$ and
$i_{h_X}(z)=-z$.  Locally around $x$ we can write $q=z^{2k}Q(z)(dz)^2$, where
$k\geq 0$ and $Q(0)\neq 0$.  Since the hyperelliptic involution $i_{h_X}$ fixes
$q$, and $i_{h_X}^*(q)=z^{2k}Q(-z)(dz)^2$, we should have that $Q(z)=Q(-z)$.
Locally around $x$, $q$ has two square roots $\pm z^k\sqrt{Q(z)}dz$.  These
square roots are fixed under $i_{h_X}(z)=-z$ if and only if $k$ is not even,
i.e. $x$ is a zero of $q$ of order not divisible by 4.

\end{proof}

The following calculation shows what happens to zeroes of abelian and
quadratic differentials when they are pulled back via a double covering.

\begin{lem}\label{lem:ZeroBehavior} Let $p:R\to S$ be a double covering, $\theta$
is an abelian differential on $S$ and $\tau$ is a quadratic differential on $S$. 

\lempart Outside of the set of branching points each zero of $\theta$ and
$\tau$ on $S$ gives rise to two zeroes of the same order on $R$.

\medskip\noindent
If $s\in S$ is a branching point of $p$, then

\smallskip\noindent
\lempart $\ord_{p^{-1}(s)}p^*(\theta)=2\ord_s\theta+1$\\
\lempart $\ord_{p^{-1}(s)}p^*(\tau)=2\ord_s\tau+2$
\end{lem}

\begin{proof}Part (a) is obvious. To prove (b) we can choose local coordinates
around $s$ and $p^{-1}(s)$ in which the map $p$ is given by $w(z)=z^2$.  If
$\theta=w^{k}S(w)dw$, $S(0)\neq 0$ and $k=\ord_s\theta$, then
$p^*(\theta)=z^{2k}S(z^2)d(z^2)=z^{2k+1}S(z^2)dz$. This proves part (b).  Part
(c) is checked similarly.
\end{proof}

\section{Veech surfaces in \texorpdfstring{$\QQM(2,2)$}{QM(2,2)}}\label{sec:QM22}

Let $(X,q)\in\QQM(2,2)$. Since $f$ is ramified only at odd-order zeroes of
$q$, the covering $f:Y\to X$ is unramified. Using Riemann-Hurwitz formula, we
obtain that $Y$ is a genus 3 surface.  \autoref{GenusY} implies that $Z$ is a
torus.  

Each of the two zeroes of $q$ are Weierstrass points. Otherwise, they would
have to be interchanged under the hyperelliptic involution, and therefore there
would exist an abelian differential with simple zeroes at both zeroes of $q$
(see \autoref{sec:HyperellipticSurfaces}).  The square of this abelian
differential would then have the same zero divisor as $q$, which would imply
that it is proportional to~$q$, contradicting our assumption that $q$ is not a
square of an abelian differential.

Applying \autoref{IgFixedPoints}, we see that each point in the fibers over
the two double zeroes of $q$ is fixed under $i_g$. The double covering $g:Y\to
Z$ is ramified at these four points. \autoref{lem:ZeroBehavior} implies that
the abelian differential $\omega$ on Z has no zeroes (which is not surprising, since
$Z$ is a torus). 

Let $B_g\subset Z$ be the 4 branching points of $g:Y\to Z$. By
\autoref{MainThm}, $SL(X,q)$ and $SL(Z,\omega,B_g)$ are commensurate. Since
$SL(Z,\omega,B_g)\subset SL(Z,\omega)\cong SL_2(\Z)$, we get the following
theorem:

\begin{thm} All Veech surfaces in $\QQM(2,2)$ are arithmetic.
\end{thm}

To find all Veech surfaces in $\QQM(2,2)$, according to \autoref{MainThm}, we
need to establish when $SL(Z,\omega,B_g)$ is a finite index subgroup
in $SL(Z,\omega)\cong SL_2(\Z)$

\begin{defn}
A finite subset $S$ in a torus $T$ is called {\it rational} if one can
identify $T$ with $\C/\Lambda$ for some lattice $\Lambda\subset\C$, so that all
points of $S$ have rational coordinates with respect to some (and therefore
any) basis of $\Lambda$.
\end{defn}

\begin{lem}\label{lem:SL(Torus,B)} 
	
	$SL(Z,\omega,B_g)$ is a finite index subgroup of $SL(Z,\omega)\cong
	SL_2(\Z)$ if and only if $B_g$ is a rational subset of $Z$.

\end{lem}
\begin{proof}

For a finite subset $F\subset Z$, denote by $Aff^+_F(Z,\omega)$ the affine
diffeomorphisms of $(Z,\omega)$ that fix $F$ pointwise, and by $SL_F(Z,\omega)$
the linear parts of these diffeomorphisms. It is clear that $SL_F(Z,\omega)$ is
a finite index subgroup of $SL(Z,\omega,F)$. Indeed, a suitable power of a
diffeomorphism permuting $F$ will fix $F$ pointwise.

Therefore $SL_{B_g}(Z,\omega)$ is a finite index subgroup of $SL(Z,\omega,B_g)$.
Hence it suffices to prove that $SL_{B_g}(Z,\omega)$ is a finite index subgroup
of $SL(Z,\omega)$ if and only if $B_g$ is a rational subset of $Z$.

Pick a point $b\in B_g$ and a lattice $\Lambda\subset\C$, so that
$Z=\C/\Lambda$ and $b\mapsto 0+\Lambda$. Under the map $\phi\mapsto D\phi$, the
group $Aff^+_{\{b\}}(Z,\omega)$ is identified with $SL(Z,\omega)\cong
SL_2(\Z)$: every orientation-preserving affine map of $(Z,\omega)$ that fixes
point $b$ can be lifted to a linear map of $\C$ preserving lattice $\Lambda$.
The group $SL_{B_g}(Z,\omega)=Aff^+_{B_g}(Z,\omega)$ is the pointwise
stabilizer of $B_g$ under the action of $SL_2(\Z)$ on $\C/\Lambda$. Therefore
it is a finite index subgroup of $SL(Z,\omega)$ if and only if every point of
$B_g$ has a finite orbit under the action of $SL_2(\Z)$, which is equivalent to
every point of $B_g$ having rational coordinates with respect to $\Lambda$. 

\end{proof}

\begin{figure}[ht]
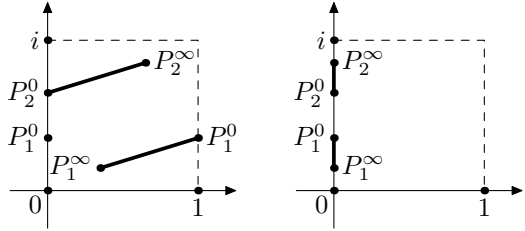

\begin{center}
$$\includegraphics{fig-3} \qquad \includegraphics{fig-4}$$
\caption{Zeroes and poles of quadratic differential $\bar q$ on the standard torus $\C/\Z^2$.}
\label{fig:TorusCuts}
\end{center}
\end{figure}

We have associated to every Veech surface in $\QQM(2,2)$ a rational four-point
subset $B_g\subset Z$. We need to reconstruct the original surface $(X,q)$ from
$(Z,\omega,B_g)$. As explained in \autoref{MainReconstruction}, this can be
achieved by choosing an abelian differential $\bar\omega$ on $X$. We can do
this so that $\bar\omega$ has a double zero at one of the zeroes of $q$ (such
$\bar\omega$ exists since both zeroes of $q$ are Weierstrass points).
\autoref{lem:ZeroBehavior} implies that the corresponding quadratic
differential $\bar q$ has simple zeroes at two points of $B_g$ and simple poles
at the other two points of $B_g$. Let $P^0_1$ and $P^0_2$ be the zeroes of
$\bar q$ and $P^\infty_1$ and $P^\infty_2$ be the poles of $\bar q$.  Note that
the involution $i_f$ descends to a hyperelliptic involution $i_{h_Z}$ of $Z$
(\autoref{ZHyperelliptic}) such that $i_{h_Z}(P^0_1)=P^0_2$ and
$i_{h_Z}(P^\infty_1)=P^\infty_2$.  Pick a lattice $\Lambda\subset\C$ so that
$Z=\C/\Lambda$ and points of $B_g$ have rational coordinates with respect to
$\Lambda$.  Hyperelliptic involution $i_{h_Z}$ of $\C/\Lambda$ lifts to
$z\mapsto-z+a$ on $\C$ for some $a\in\C$. Since $i_{h_Z}$ permutes $B_g$, $a$
has rational coordinates with respect to $\Lambda$. By composing projection
$\pi:\C\to\C/\Lambda$ with a translation on the torus, we can assume that $a=0$
and $B_g$ still has rational coordinates with respect to $\Lambda$. Since
$i_{h_Z}(P^0_1)=P^0_2$, three points $P^0_1$, $P^0_2$ and $\pi(0)$ are
colinear. Moreover, since coordinates of these three points are rational with
respect to $\Lambda$, any geodesic through these three points is closed. Since
we are interested in describing Veech surfaces up to the action of $SL_2(\R)$,
we can act by a suitable element from $SL_2(\R)$ to map $\Lambda$ into the
standard lattice $\Z^2$, so that the line through the three points  $P^0_1$,
$P^0_2$ and $\pi(0)$ is the image of the imaginary axis. Depending on whether
all four points $P^0_1, P^0_2, P^\infty_1, P^\infty_2$ are colinear, we will
get two possible configurations illustrated in \autoref{fig:TorusCuts}.

Zeroes of $\bar q$, as well as poles of $\bar q$, are symmetric about all 4 points of order 2 on the torus $\C/\Z^2$.

Now assume we start with a standard torus $Z=\C/\Z^2$ and four rational points
$P^0_1$, $P^0_2$, $P^\infty_1$ and $P^\infty_2$ as above. There exists unique
(up to multiplication by a scalar) meromorphic differential $\bar q$ having
simple zeroes at $P^0_1$, $P^0_2$ and simple poles at $P^\infty_1$,
$P^\infty_2$. Indeed, involution $i_{h_Z}:z\mapsto -z$ defines a 2-to-1 map
$h_Z:Z\to\CP^1$. By composing it with a Moebius transformation on $\CP^1$ we
can assume that $h_Z(P^0_1)=h_Z(P^0_2)=0$ and
$h_Z(P^\infty_1)=h_Z(P^\infty_2)=\infty$. Then $\bar q=h_Z(z)\,dz^2$. It is
clear that $\bar q$ is fixed under the hyperelliptic involution $z\mapsto-z$.
Hence we can apply the main construction to $(Z,\bar q)$ in order to describe
$(X,q)$.

\begin{figure}[ht]
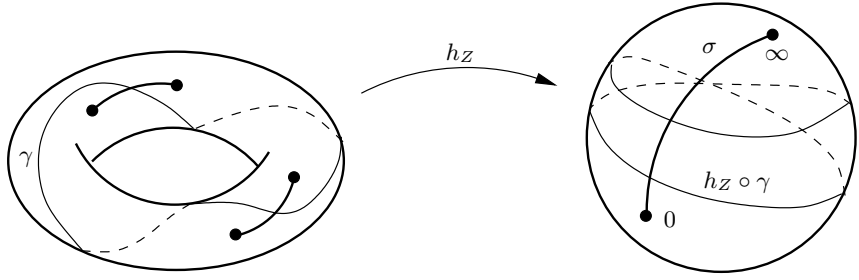

\begin{center}

\ifx\pdftexversion\undefined
\input{pics/pic1.pstex_t}
\else
\input{pics/pic1.pdftex_t}
\fi

\caption{Curve $\gamma$ with trivial monodromy.}
\label{fig:TorusMonodromy}
\end{center}
\end{figure}

To describe the covering $g:Y\to Z$ we need to see when the monodromy along a
closed curve $\gamma:[0,1]\to Z$ avoiding branching points $P^0_1, P^0_2,
P^\infty_1, P^\infty_2$ is trivial. Since $g:Y\to Z$ is the double covering
given by $\bar q=h_Z(z)\,dz^2$, the monodromy along $\gamma$ will be trivial if
the curve $h_Z\circ\gamma:[0,1]\to\C$ has even winding number around zero.
If we cut $\CP^1$ from zero to infinity along some simple path $\sigma$ avoiding
branching points of $h_Z$, then the monodromy of $\gamma$ is trivial if and
only if $h_Z\circ\gamma$ intersects $\sigma$ even number of times or,
equivalently, if $\gamma$ intersects $h_Z^{-1}(\gamma)$ even number of times (see \autoref{fig:TorusMonodromy}).

We can make parallel cuts on $Z$ as shown in \autoref{fig:TorusCuts}. Under
$h_Z$ these cuts get identified into one curve connecting zero to infinity.

From the discussion above we see that we can think of $Y$ as two copies of the
torus $Z$ glued along the cuts, so that when we cross one of the cuts we go to
the other copy of $Z$: $Y=(Z_1\sqcup Z_2)/(\partial Z_1\sim\partial Z_2)$.

Now that we understand the structure of $Y$, we can explain how to get $X$
from~$Z$. $X$ is obtained from $Y$ by factoring by the involution $i_f:Y\to Y$.
By \autoref{ZHyperelliptic}, $i_f$ descends to the hyperelliptic involution
$i_{h_Z}$ on $Z$. We have also seen that $i_f$ has no fixed points. This
implies that $i_f$ has to interchange $Z_1\backslash\partial Z_1$ and
$Z_2\backslash\partial Z_2$. Indeed, say a point $A$ from
$Z_1\backslash\partial Z_1$ was mapped to another point $B\in
Z_1\backslash\partial Z_1$ symmetrical about the center of the unit square.
Then we can connect $A$ to $B$ avoiding the cuts by a path symmetrical about
the center of the unit square (the path does not have to be straight). By
continuity, along this path $i_f$ will have to equal $i_{h_{Z_1}}$, which will
mean that the center of the square is fixed under $i_f$, contradicting that
$i_f$ has no fixed points. Since $i_f$ glues interiors of two copies of $Z$, we
just need to see how it acts on the sides of the cuts. As easily verified, we
get the following statement:

\begin{thm} Up to the action of $SL_2(\R)$ any Veech surface in $\QQM(2,2)$ can be obtained by making two non-intersecting (but possibly colinear) parallel cuts of the same length with rationally related vertices on the standard torus $(\C/\Z^2,dz)$ and gluing as indicated in \autoref{fig:QM22}.
	
\begin{figure}[ht]
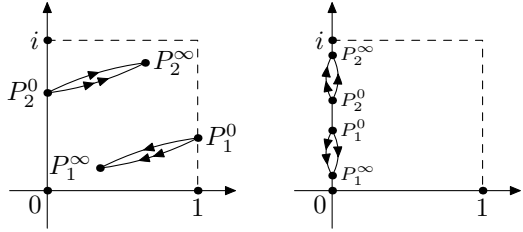

\begin{center}
$$\includegraphics{fig-5} \qquad \includegraphics{fig-6}$$
\caption{Veech surface in $\QQM(2,2)$. Non-colinear (left) and colinear (right) cuts.}
\label{fig:QM22}
\end{center}
\end{figure}
\end{thm}

It is easy to check that there will be infinitely many configurations
of four points $\{P^0_1, P^0_2, P^\infty_1, P^\infty_2\}$ up to the action of
$SL_2(\R)$ on the standard torus $\C/\Z^2$. Therefore there are infinitely many
Veech surfaces in $\QQM(2,2)$.

\section{Veech surfaces in \texorpdfstring{$\QM(2,1,1)$}{QM(2,1,1)}}

Assume quadratic differential $q$ has one double zero and two
simple zeroes. Then the covering $f:Y\to X$ is branched over two simple zeroes
of $q$.  Riemann-Hurwitz formula implies that $Y$ is a genus 4 surface, and
therefore by \autoref{GenusY} $Z$ is a genus 2 surface.

\begin{figure}[ht]
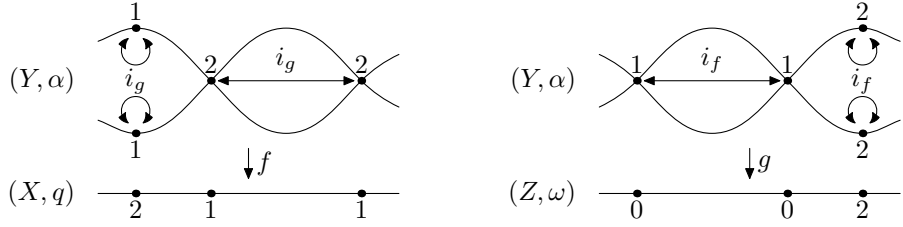

$$\includegraphics{fig-7} \qquad\qquad \includegraphics{fig-8}$$

\caption{Structure of the covering maps $f:Y\to X$ and $g:Y\to Z$, and
corresponding involutions $i_f$ and $i_g$, along with the orders of zeroes of
the forms $q$, $\alpha$ and~$\omega$.}

\label{fig:QM211CoveringStructure}

\end{figure}

Hyperelliptic involution $h_X$ preserves $q$, hence the double zero of $q$ is a
Weierstrass point. The simple zeroes of $q$ are conjugate under the
hyperelliptic involution. Indeed, if they were fixed under the hyperelliptic
involution, then they would have to be Weierstrass points. Since $q$ is a
product of two abelian differentials, the order of zero of $q$ at a Weierstrass
point has to be even (see \autoref{sec:HyperellipticSurfaces}), therefore $q$
can not have a simple zero at a Weierstrass point.

By \autoref{IgFixedPoints}, $i_g$ has two fixed points on $Y$, namely the
points in the fiber over the double zero of $q$. Applying
\autoref{lem:ZeroBehavior}, we get that abelian differential $\alpha$ on $Y$
has two double zeroes over the simple zeroes of $q$, and two simple zeroes over
the double zero of $q$. Consequently, abelian differential $\omega$ on $Z$ has
one double zero.

Let $B_g$ be the set of branching points of $g:Y\to Z$. As indicated on the
figure above $B_g$ consists of two points. According to \autoref{MainThm},
$(X,q)$ is Veech if and only if $(Z,\omega)$ is Veech and $SL(Z,\omega,B_g)$ is
a finite index subgroup of $SL(Z,\omega)$. Assume $(X,q)$ and $(Z,\omega)$ are
Veech. The following lemma shows that the two points of $B_g$ have to be
periodic under the action of $Aff^+(Z,\omega)$.

\begin{lem}\label{lem:SL(H,B)} Let $\tau$ be a quadratic differential on a
	closed Riemann surface $R$ of genus at least 2, and $B$ be a finite
	subset of $R$. Then $SL(R,\tau,B)$ is a finite index subgroup of
	$SL(R,\tau)$ if and only if $B$ has a finite orbit under the action of
	$Aff^+(R,\tau)$ or equivalently, that every point of $B$ is periodic
	under the action of $Aff^+(R,\tau)$. 

\end{lem}

\begin{proof} 

Consider the exact sequence: $$0\lto Aut(R,\tau)\lto
Aff^+(R,\tau)\stackrel{D}{\lto}SL(R,\tau)\to 0 $$ Since $R$ is a genus 2
surface, there are finitely many holomorphic automorphisms of $R$, and
therefore $Aut(R,\tau)$ is finite.  The group $Aff^+(R,\tau,B)$ is mapped onto
$SL(R,\tau,B)$ under $D$. Using the exact sequence above and that
$Aut(R,\tau)$ is finite, we see that $SL(R,\tau,B)$ is a finite index subgroup
of $SL(R,\tau)$ if and only if $Aff^+(R,\tau,B)$ is a finite index subgroup of
$Aff^+(R,\tau)$. Since $B$ is a finite set, $Aff^+_B(R,\tau)$, the pointwise
stabilizer of B under the action of $Aff^+(R,\tau)$, is a finite index subgroup
of $Aff^+(R,\tau,B)$. The group $Aff^+_B(R,\tau)$ has finite index in
$Aff^+(R,\tau)$ if and only if $B$ has a finite orbit under the action of the
affine group $Aff^+(R,\tau)$.
\end{proof}

Periodic points on translation surfaces have been studied in a paper by Gutkin,
Hubert and Schmidt (\cite{GHS}). They prove that under certain conditions all
Weierstrass points on a hyperelliptic surface are periodic. In a more recent
paper, M\"oller showed that the converse is always true in the case of a
primitive Veech surface of genus 2 in $\Omega M_2(2)$.

\begin{thm}[M\"oller, {\cite[Theorem 5.1]{Mo}}]\label{Moller}
The only periodic points on a primitive Veech surface in $\Omega M_2(2)$ are the Weierstrass points.
\end{thm}

Since $i_f$ descends to the hyperelliptic involution $h_Z$ on $Z$, two points
of $B_g$ are interchanged under $h_Z$. Hence they are not Weierstrass. This
along with the \autoref{Moller} means that $(Z,\omega)$ cannot be
primitive, thus it has to be a cover of a torus. Since $SL(X,q)$ and $SL(Z,\omega)$ are commensurate, we get the following theorem: 

\begin{thm}All Veech surfaces in $\QM(2,1,1)$ are arithmetic.
\end{thm}

\section{Veech surfaces in \texorpdfstring{$\QM(1,1,1,1)$}{QM(1,1,1,1)}}

\subsection{General situation}
Assume quadratic differential $q$ has four simple zeroes on $X$. Then the map
$f:Y\to X$ is branched over four points. Therefore $Y$ has genus 5. The map
$g:Y\to Z$ has no branching points, and hence $Z$ has genus 3. The holomorphic
form $\alpha$ has four zeroes of order 2 on $Y$, which are mapped under $g:Y\to
Z$ into two double zeroes of $\omega$. Both zeroes of $\omega$ are Weierstrass
points of $Z$.

Since the covering $g:Y\to Z$ is not ramified, \autoref{MainThm} implies
that $(X,q)$ is Veech if and only if $(Z,\omega)$ is Veech. This way we get a
map from $\QM(1,1,1,1)$ to $\OM[3](2,2)$, which sends Veech surfaces to
hyperelliptic Veech surfaces that have singularities at Weierstrass points.

\begin{thm}\label{thm:QM1111}
There is a one-to-one correspondence between Veech surfaces of genus 2 in $\QM(1,1,1,1)$ and hyperelliptic Veech surfaces of genus 3 in $\OM[3](2,2)$ with singularities at Weierstrass points.
\end{thm}
\begin{proof}

Every quadratic differential on a genus 2 surface is a product of abelian
differentials (see \autoref{sec:HyperellipticSurfaces}). Therefore the quadratic
differential $q$ is a product of two abelian differentials $\bar\omega_1$ and
$\bar\omega_2$ of type $(1,1)$. Take one of them, say $\bar\omega_1$ and apply
to it the construction from \autoref{MainReconstruction}. We will get a
quadratic differential $\bar q$ on $Z$ with a zero of order 6 at one of the
zeroes of $\omega$ and a zero of order 2 at the other zero of $\omega$. If we
choose $\bar\omega_2$ instead of $\bar\omega_1$ at the beginning, then the
zeroes of $\bar q$ will be switched.

Now suppose we start with any hyperelliptic genus 3 surface $Z$ with an abelian
differential $\omega$ of type $(2,2)$, where both zeroes are Weierstrass
points. There exists an abelian differential that has a zero of order 4 at one
of the zeroes of $\omega$. Let $\bar q$ be the product of this differential and
$\omega$. Then $\bar q$ is a quadratic differential of type $(6,2)$, and it is
the only quadratic differential with such zero divisor up to multiplication by
a scalar. Moreover $\bar q$ is not a square of an abelian differential, because
an abelian differential on a hyperelliptic surface cannot have a zero of odd
order at a Weierstrass point.  Furthermore, $\bar q$ is fixed by the
hyperelliptic involution because it is a product of abelian differentials.
Applying ideas from \autoref{MainReconstruction} to $(Z,\omega,\bar q)$ we get
a triple $(X,q,\bar\omega)$. It can be easily verified using results from
\autoref{MainConstruction}, that $X$ is a genus 2 surface and $q$ is a
quadratic differential with four simple zeroes.

\end{proof}

We will use the following description of a double covering given by a quadratic
differential to better understand the correspondence between Veech surfaces
in $\QM(1,1,1,1)$ and $\OM[3](2,2)$.

\begin{lem}\label{lem:HyperellipticQDCuts}

Assume $X$ is a hyperelliptic Riemann surface of genus at least~2 with the
hyperelliptic involution $i:X\to X$. Let $q$ be a quadratic differential on $X$
with two zeroes of even orders at Weierstrass points $P_1$ and $P_2$. Assume
additionally that $q$ is not a square of an abelian differential. Let $Y\to X$
be the double covering given by $q$.  Connect $P_1$ to $P_2$ by a simple path
$\sigma$ not passing through other Weierstrass points of $X$. Then $Y$ can be
obtained by taking two copies of $X$, cutting them along a closed loop
$\sigma\cup i(\sigma)$ and re-gluing them in the standard way to obtain a
double cover of $X$. 

\end{lem}

\begin{proof} 
	
Let $P_1,P_2,\ldots,P_{2g+2}$ be all Weierstrass points of $X$. Let
$z:X\to\CP^1$ be the double cover defined by the hyperelliptic involution $i$,
chosen so that $z(P_j)\neq\infty$, $j=1,\ldots,2g+2$. Consider the function
$w=\sqrt{\prod\limits_{j=1}^{2g+2}(z-z(P_j))}$. Then $X$ is the Riemann surface
defined in complex coordinates $(z,w)$ by
$w^2=\prod\limits_{j=1}^{2g+2}(z-z(P_j))$. 

Assume $deg_{P_1}(q)=2a$ and $deg_{P_2}(q)=2b$, $a,b\in\Z$. The total number of
zeroes of $q$ is $2a+2b=4g-4$. Therefore $a+b$ is even, and $a$ and $b$ have
the same parity. If they were both even, then $q$ would be a square of an
abelian differential with the zero divisor $aP_1+bP_2$. Hence both $a$ and $b$
are odd. 

In coordinates $(z,w)$ we can express

\begin{equation}
\label{eq:q}q=\frac{(z-z(P_1))^a(z-z(P_2))^b(dz)^2}{w^2}=\frac{(z-z(P_1))^{a-1}(z-z(P_2))^{b-1}(dz)^2}{(z-z(P_3))\ldots(z-z(P_{2g+2}))}
\end{equation}

Take any closed path $\gamma$ on $X$. We would like to see when the monodromy
of the covering $Y\to X$ along $\gamma$ is non-trivial, or equivalently, when
the extension of $\sqrt{q}$ along $\gamma$ changes sign. Since $a$ and $b$ are
odd, the numerator in (\ref{eq:q}) is a square of a well-defined abelian
differential.  Therefore the monodromy along $\gamma$ is non-trivial if and
only if the analytical extension of $\sqrt{(z-z(P_3))\ldots(z-z(P_{2g+2}))}$
along $\gamma$ changes sign. This will happen if and only if the total winding
number of $z(\gamma)$ about points $z(P_3),\ldots,z(P_{2g+2})$ is odd. The
winding number of $z(\gamma)$ about the images of all Weierstrass points
$z(P_1),\ldots,z(P_{2g+2})$ is always even (this is the reason why the function
$w:X\to\CP^1$ is well-defined). Therefore the winding number of $z(\gamma)$
about $z(P_3),\ldots,z(P_{2g+2})$ is odd if and only if the winding number of
$z(\gamma)$ about $z(P_1)$ and $z(P_2)$ is odd.

We are given a path $\sigma$ connecting $P_1$ to $P_2$. The winding number of
$z(\gamma)$ about $z(P_1)$ and $z(P_2)$ is odd if and only if the intersection
number of $z(\gamma)$ and $z(\sigma)$ is odd, or equivalently the intersection
number of $\gamma$ and $\sigma\cup i(\sigma)$ is odd. This finishes the proof
of the lemma.

\end{proof}

Consider a genus 3 surface $(Z,\omega)$ and a quadratic differential $\bar q$
as in \autoref{thm:QM1111}. Denote zeroes of $\omega$ (and $\bar q$) by
$P_1$ and $P_2$. Connect $P_1$ to $P_2$ by a saddle connection $\sigma$.
Following the construction in \autoref{thm:QM1111}, we get a genus 5
surface $Y$ which is a double cover of $Z$ given by $\bar q$. The pair $(Z,\bar
q)$ satisfies the hypothesis of \autoref{lem:HyperellipticQDCuts}. Therefore
$Y$ is obtained by taking two copies of $Z$, cutting them along $\sigma\cup
i_{h_Z}(\sigma)$ and gluing together along the cuts: $Y=(Z_1\sqcup
Z_2)/(\partial Z_1\sim\partial Z_2)$, where $\partial Z_i=\sigma_i\cup
i_{h_{Z_i}}(\sigma_i),\ i=1,2$.

Surface $(X,q)$ is obtained by factoring $Y$ by the involution $i_f$. By
\autoref{ZHyperelliptic}, $i_f$ descends to the hyperelliptic involution
$i_{h_Z}$ on $Z$. By \autoref{IgFixedPoints}, $i_f$ interchanges the
Weierstrass points of $Z_1$ and $Z_2$ that are not zeroes of $\bar q$. Surfaces
$Z_1$ and $Z_2$ are connected, because $Y$ is (equivalently, the curve
$\sigma\cup i_{h_Z}(\sigma)$ is non-separating).  Therefore by continuity, $i_f$
interchanges interiors of $Z_1$ and $Z_2$. On the sides of the cuts $\partial
Z_i$, $i_f$ is given by the hyperelliptic involution $i_{h_{Z_i}}$. We have
proved the following statement:

\begin{figure}[ht]
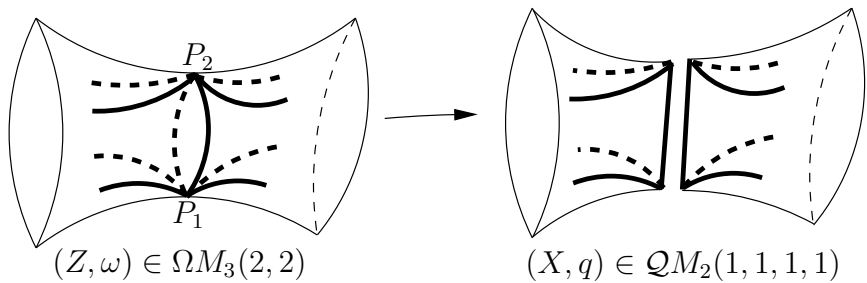

\begin{center}

\ifx\pdftexversion\undefined
\input{pics/pic2.pstex_t}
\else
\input{pics/pic2.pdftex_t}
\fi

\caption{Correspondence between Veech surfaces in $\OM[3](2,2)$ and $\QM(1,1,1,1)$. Bold lines indicate the foliation in the direction of a saddle connection from $P_1$ to $P_2$.}
\label{fig:OM22toQM1111}
\end{center}
\end{figure}
\begin{thm}\label{thm:QM1111v2}
	
Every Veech surface in $\QM(1,1,1,1)$ can be obtained from a hyperelliptic
Veech surface in $\OM[3](2,2)$ with singularities at the Weierstrass points in
the following way: take a saddle connection joining the two zeroes, act on it
with the hyperelliptic involution to obtain a closed loop, then cut the surface
along this loop and glue each opening shut via the hyperelliptic involution
(see \autoref{fig:OM22toQM1111}).

\end{thm}

\subsection{Examples}

While the author is not aware of any primitive genus 3 Veech surfaces
satisfying conditions of \autoref{thm:QM1111}, we can construct some
non-primitive, non-arithmetic examples by considering an unramified double
covering over a primitive genus~2 Veech surface in $\OM(2)$. 

\begin{figure}[ht]
\begin{center}
\includegraphics{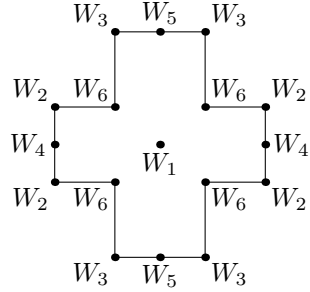}

\caption{Weierstrass points on a Swiss cross surface $(S,\theta)\in\OM(2)$ arising from an L-shaped billiard table. $W_6$ is the double zero of $\theta$.}\label{fig:SwissCrossWP}

\end{center}
\end{figure}

According to
\cite{McMSpin}, all such surfaces arise from L-shaped billiard tables of
specific dimensions (see \autoref{IntroVeech2}).  The surface itself is
obtained by reflecting an L-shaped table 4 times across its sides until we get a
'Swiss cross', and then gluing the parallel sides. The hyperelliptic involution is given by rotating by $\pi$ about the center of the cross. The Weierstrass points are shown on \autoref{fig:SwissCrossWP}. 

Choose a Veech surface $(S,\theta)\in\OM(2)$. An unramified double covering
$(Z,\omega)\to (S,\theta)$ corresponds to a choice of a quadratic differential
$\tau$ on $S$ with two double zeroes.  We are interested in non-trivial double
coverings, therefore $\tau$ should not be a square of an abelian differential.
This implies that the two zeroes of $\tau$ have to occur at the Weierstrass
points of $S$ (see the beginning of \autoref{sec:QM22}). To satisfy the
hypothesis of \autoref{thm:QM1111}, we need to have that the abelian
differential $\omega$ has its zeroes at the Weierstrass points of $Z$. By
\autoref{cor:YHyperelliptic} and \autoref{IgFixedPoints}, Weierstrass points of
$Z$ are the eight points lying over those Weierstrass points of $S$ that are
not zeroes of the quadratic differential~$\tau$. Therefore $\tau$ cannot have
one of its zeroes at $W_6$, which is the zero of $\theta$. Hence zeroes of
$\tau$ can only occur at the five Weierstrass points $W_1,\ldots,W_5$. Having
chosen two such Weierstrass points on $S$, we use
\autoref{lem:HyperellipticQDCuts} to obtain a concrete description of the
surface $(Z,\omega)\in\OM[3](2,2)$. Then we use \autoref{thm:QM1111v2} to
obtain a Veech surface $(X,q)\in\QM(1,1,1,1)$.

There are ten different ways to choose two zeroes of $\tau$, hence there are
at most ten different Veech surfaces in $\QM(1,1,1,1)$ corresponding to a given
Veech surface in $\OM(1,1)$. 

We will look in detail at the case in which $W_1$ and $W_2$ are the zeroes of
$\tau$.

\begin{figure}[ht]
\begin{center}
\includegraphics{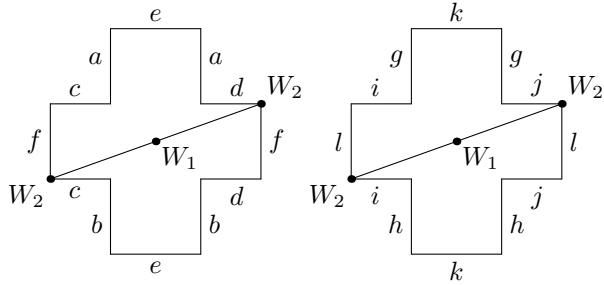}

\caption{Two copies of a Veech surface $(S,\theta)\in\OM(2)$. Lowercase letters
indicate which sides are identified.}\label{fig:SwissCross1}

\end{center}
\end{figure}
	
\begin{figure}[ht]
\begin{center}
\includegraphics{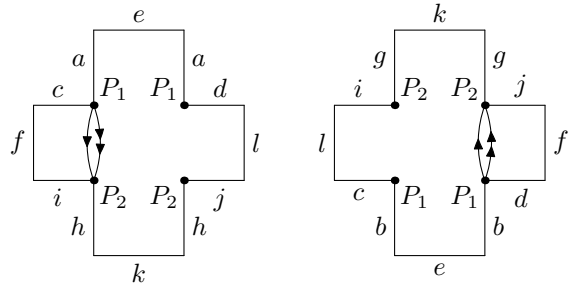}

\caption{Genus 3 surface $(Z,\omega)\in\OM[3](2,2)$.}\label{fig:SwissCross2}

\end{center}
\end{figure}

Take two copies of a Veech surface $(S,\theta)\in\OM(2)$ (see
\autoref{fig:SwissCross1}). Following \autoref{lem:HyperellipticQDCuts}, pick a
path from $W_1$ to $W_2$ and act on it with the hyperelliptic involution to
obtain a closed loop on each surface. Cut the surfaces along these loops and
re-glue the surfaces along the loops switching the order to get a non-trivial
double cover of $S$. This results in a genus 3 surface displayed in
\autoref{fig:SwissCross2}.

\begin{figure}[ht] \begin{center} \includegraphics{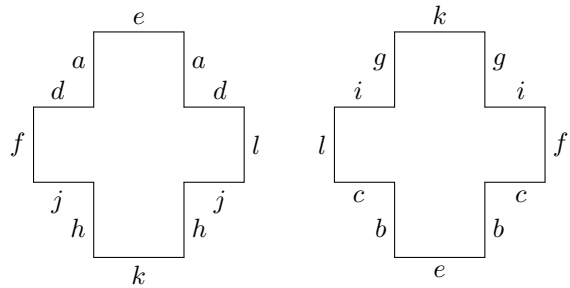}

\caption{Genus 2 Veech surface $(X,q)\in\QM(1,1,1,1)$.}\label{fig:SwissCross3}

\end{center}
\end{figure}

Points $P_1$ and $P_2$ in \autoref{fig:SwissCross2} are the two zeroes of
$\omega$. By \autoref{cor:YHyperelliptic}, the hyperelliptic involution on
$(Z,\omega)$ projects to the hyperelliptic involution on $(S,\theta)$. As
mentioned above, Weierstrass points of $Z$ are the eight points lying over the
four Weierstrass points $W_3,W_4,W_5,W_6$.  Hence the hyperelliptic involution
of $Z$ has no fixed points in the interiors of the two Swiss crosses. Therefore
the hyperelliptic involution acts on $Z$ by rotating each Swiss cross about its
center by $\pi$ and then switching the sheets. Following
\autoref{thm:QM1111v2}, pick a path from $P_1$ to $P_2$ and act on it with the
hyperelliptic involution to obtain a closed loop.  The two sides of this loop
are indicated in \autoref{fig:SwissCross2} using single and double arrows.  We
need to cut the surface along this loop and glue each opening shut via the
hyperelliptic involution. This corresponds to identifying the segments marked
by single arrows and the segments marked by double arrows.  The resulting
surface is displayed in \autoref{fig:SwissCross3}. 

\begin{figure}[ht]
\begin{center}
\includegraphics{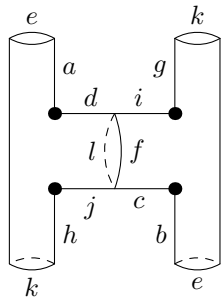}

\caption{Genus 2 Veech surface $(X,q)\in\QM(1,1,1,1)$. Sides marked with $e$
and sides marked with $k$ are identified with a twist by $\pi$ so that the
segments $a$ and $b$, and $g$ and $h$ are aligned. The four marked points are
the zeroes of~$q$.}\label{fig:SwissCross4}

\end{center}
\end{figure}
After performing some of the identifications in \autoref{fig:SwissCross3}, we
will get an H-shaped surface displayed in \autoref{fig:SwissCross4}. Provided
that we started with a non-arithmetic Veech surface in $\OM(2)$,
\autoref{fig:SwissCross4} is an example of a non-arithmetic Veech surface in
$\QM(1,1,1,1)$.

\FloatBarrier

\end{document}